\documentclass[12pt]{article}      
\usepackage{amsmath}
\usepackage{amsfonts}
\usepackage{amsthm}
\usepackage{amssymb}
\usepackage{mathtools}
\usepackage[margin=1in]{geometry}
\usepackage{graphicx}
\usepackage{graphics}
\usepackage{listings}
\usepackage{multirow, tabu}
\usepackage{framed}
\usepackage{hyperref}
\usepackage{amsfonts}
\usepackage{amstext}
\usepackage{array}
\usepackage{float} 
\usepackage{verbatim}
\usepackage{subfigure}
\usepackage{tikz-cd}
\usepackage{tikz}
\usetikzlibrary{arrows,automata,shapes,graphs,fit,patterns} 
\usetikzlibrary{positioning}
\usetikzlibrary{calc}
\usetikzlibrary{decorations.pathreplacing}

\newcolumntype{C}{>{$}c<{$}}

\newcommand{\F}{\mathbb{F}}
\newcommand{\N}{\mathbb{N}}
\newcommand{\Z}{\mathbb{Z}}

\usepackage{color}

\usepackage[algo2e,ruled,vlined,]{algorithm2e}
\usepackage{algorithmic}

\definecolor{red}{rgb}{1,0,0}

\definecolor{blue}{rgb}{0,0,1}

\usepackage{mathrsfs}

\newcommand{\ms}{\medskip}
\newcommand{\bpf}{\begin{proof}}
\newcommand{\epf}{\end{proof}\ms}

\usepackage{lineno}

\newtheorem{theorem}{Theorem}

\newtheorem{lemma}[theorem]{Lemma}

\newtheorem{observation}[theorem]{Observation}

\theoremstyle{definition}
\newtheorem{definition}{Definition}

\tikzset{addbrace/.style = {decorate,decoration={brace,raise=5pt, amplitude=1.0ex,}, thick}}

\title{Templates, Arrays, and Overlays}
\date{}
\author{Jordan Broussard}
\begin{document}
\maketitle


An \emph{array} is a function $A: \Z^2 \rightarrow \F$, where $\F$ is a field.  We adopt the matrix notation $A_{r,c}$ instead of $A(r,c)$.  An array can be thought of as an extention of Riordan arrays (removing the condition that the array is lower-triangular) with indices from $\Z$ rather than indices from $\N$.  Instead of using the rows of the array to store powers of the coefficients of the generating function\cite{Wilf}, we store the coefficients of a generating function in two variables ($x$ and $y$) where the degree of $x$ denotes the row and the degree of $y$ denotes the column. ***Add comment about how to actually use it as a Riordan array *** Since the array allows negative indices, then the corresponding generating function is a formal Laurent series rather than a formal power series.  That is, we can define the generating function $a(x,y)$ as $$a(x,y)=\sum_{(r,c) \text{ } \in \text{ } \Z^2} A_{r,c} \text{ }x^r y^c.$$ \cite{FlajSedg}


\begin{observation} The set of arrays is a vector space under entry-wise addition and scalar multiplication. \end{observation}

Let $X$ and $Y$ be shift operators such that $(XA)_{r,c}=A_{r,c-1}$ and $(YA)_{r,c} = A_{r-1,c}$.  Notice that $X$ and $Y$ are invertible linear operators and they commute.  

\begin{definition}
A \emph{template $T$} is a polynomial in $X$ and $Y$, where $X^0=Y^0=I$, the identity operator.
\end{definition}

Since $T$ is a polynomial of linear operators, then $T$ is also linear.  

\begin{definition}
$A$ is \emph{annihilated by} the operator $\Psi$ if $\Psi A=0$, where 0 denotes the zero array.
\end{definition}

\begin{definition}
We also say that $A$ is a \emph{template array with template $T$} if $A$ is annihilated by a template $T$. 
\end{definition} 

As an example, let $A$ be a template array with template $T$, where $T=I-Y(I+X)$.  Then we can recover the two-dimensional recurrence by considering the $(r,c)$ entry of $T$ applied to $A$.  Since $A$ is a template array, then we have
\begin{align*}
(TA)_{r,c}&=0 \Leftrightarrow \\
((I-Y(I+X))A)_{r,c}&=0 \Leftrightarrow \\
((I-Y-XY)A)_{r,c}&=0 \Leftrightarrow \\
A_{r,c}-A_{r-1,c}-A_{r-1,c-1}&=0 \Leftrightarrow \\
A_{r,c}=A_{r-1,c}+A_{r-1,c-1}&.
\end{align*}

\begin{definition} 
An \emph{overlay for a template T} is an array with finite support that contains the coefficients of $T$.
\end{definition}

We typically draw overlays with only their nonzero entries visible.  We can define the \emph{overlay product} of two arrays by the Hadamard product (or the sum of entry-wise products).  We define the dimensions of an overlay to be (number of rows containing non-zero entries) $\times$ (number of columns containing non-zero entries).  Let the array $B=\left(b_{i,j}\right)$ be an overlay for the template $T$.  Let the right-most nonzero entry have column index 0 and the left-most nonzero entry have column index $n$.  Likewise, let the vertically lowest nonzero entry have row index 0 and the verically highest entry have row index $m$.  The coefficient $b$ for the term $bX^cY^r$ is the entry in the cell with row index $r$ and column index $c$.

For example, for the template $T=I-Y(I+X)=I-Y(I+X)$, the corresponding overlay would be
\begin{center}
\begin{tikzpicture}
\draw [black, very thin] (0,1.5)--(1.5, 1.5)--(1.5,0)--(.75,0)--(.75,.75)--(0,.75)--cycle;
\draw [black, very thin] (.75,1.5)--(.75,.75)--(1.5, .75);
\draw [thin, black] node at (.375, 1.125) {$-1$};
\draw [thin, black] node at (1.125, 1.125) {$-1$};
\draw [thin, black] node at (1.125, .375) {$1$};
\end{tikzpicture}
\end{center}

Let $A$ be an array annihilated by $T$.  Then the Hadamard product of the overlay $B$ (derived from the template $T$) and the array $A$ will always be zero for all $(r,c)$ such that $A_{r,c}$ is lined up to mutliply by $b_{0,0}$.  

Let $B=\left(b_{i,j}\right)$ be an overlay with dimension $(m+1) \times (n+1)$ for a template $T$ that annihilates $A$.  Then for all $r,c \in \Z$,

$$\sum_{i=0}^m \sum_{j=0}^n b_{i,j}A_{r-i,c-j} = 0.$$

For an overlay $B=\left(b_{i,j}\right)$, we define $n,m,u,$ and $l$ such that the top row of the overlay is $u+1$ units wide, the bottom row of the overlay is $l+1$ units wide, the furthest column to the right containing a nonzero entry in the overlay is indexed zero, the furthest column to the left containing a nonzero entry in the overlay is indexed $n$, and the height of the overlay is $m+1$ units tall.  We often omit writing the zero entries prior to the first nonzero entry in that row and after the last nonzero entry in that row.  An overlay may look something like the following, where $*$ is an element from $\F$.

\begin{center}
\begin{tikzpicture}
\draw [step=.75, black, very thin] (0,.75) grid (11.25,3);
\draw [step=.75, black, very thin] (2.25,3) grid (4.5, 3.75);
\draw [step=.75, black, very thin] (6.749,0) grid (9, .75);

\foreach \i in {0,1.5,2.25,3.75,4.5,6,6.75,8.25,9,10.5} {
	\foreach \j in {.75,2.25} {
		\draw [thin,black] node at (\i+.375, \j+.375) {$*$};
	}
	\draw [thin,black] node at (\i+.375,1.95) {$\vdots$};
}

\foreach \i in {.75,3,5.25,7.5,9.75} {
	\foreach \j in {.75,2.25} {
		\draw [thin,black] node at (\i+.375, \j+.375) {$\cdots$};
	}
	\draw [thin, black] node at (\i+.375, 1.95) {$\ddots$};
}

\draw [thin, black] node at (2.625,3.375) {$*$};
\draw [thin, black] node at (4.125,3.375) {$*$};
\draw [thin, black] node at (3.375,3.375) {$\cdots$};

\draw [thin, black] node at (7.125, 0.375) {$*$};
\draw [thin, black] node at (8.625, 0.375) {$*$};
\draw [thin, black] node at (7.875, 0.375) {$\cdots$};

\end{tikzpicture}
\end{center}
%
%

For example, consider the recurrence $A_{r-1,c-1}+3A_{r-1,c}+2A_{r,c-1}-A_{r,c}=0$.  Then the corresponding template would be $T=XY+3Y+2X-I$, and the corresponding overlay is

\begin{center}

\begin{tikzpicture}
\draw [step=.75, black, very thin] (0,0) grid (1.5,1.5);

\draw [thin, black] node at (.375,.375) {$2$};
\draw [thin, black] node at (1.125,.375) {$-1$};

\draw [thin, black] node at (.375,1.125) {$1$};
\draw [thin, black] node at (1.125, 1.125) {$3$};

\end{tikzpicture}
\end{center}

Now suppose that we want to construct an array $A$ that is annihilated by $T$, and we are given that $A_{0,0}=1, A_{0,c} = 0$ if $c \ne 0, A_{r,0} = 0$ if $r \ne 0$.  Then initially, $A$ looks like

\begin{center}
\begin{tikzpicture}
\draw [step=.75, black, very thin] (0,0) grid (5.25,5.25);

\draw [thin, black] node at (.375,3.375) {$\cdots$};
\draw [thin, black] node at (1.125, 3.375) {$0$};
\draw [thin, black] node at (1.875, 3.375) {$1$};
\draw [thin, black] node at (2.625, 3.375) {$0$};
\draw [thin, black] node at (3.375, 3.375) {$0$};
\draw [thin, black] node at (4.125, 3.375) {$0$};
\draw [thin, black] node at (4.875, 3.375) {$\cdots$};

\draw [thin, black] node at (1.875, .375) {$\vdots$};
\draw [thin, black] node at (1.875, 1.125) {$0$};
\draw [thin, black] node at (1.875, 1.875) {$0$};
\draw [thin, black] node at (1.875, 2.625) {$0$};
\draw [thin, black] node at (1.875, 4.125) {$0$};
\draw [thin, black] node at (1.875, 4.875) {$\vdots$};

\end{tikzpicture}

%
%
\end{center}

Imagine placing the overlay on top of $A$ such that the one in the upper left corner of the overlay is in the same position as the 1 in $A$. We can observe that three of the four corresponding entries in $A$ are known.  Since the Hadamard product must be zero, we can solve for the value in $A$ that corresponds to the bottom right corner of where the overlay is placed.  Call this value $z$.  Then by solving the equation $1*1+0*3+0*2+z*-1=0$ for $z$, we have that $z=1$.  By placing it into $A$, that gives us

\begin{figure}[H]
\centering
	\subfigure{
		\begin{tikzpicture}
		\draw [step=.75, black, very thin] (0,0) grid (5.25,5.25);

		\draw [thin, black] node at (.375,3.375) {$\cdots$};
		\draw [thin, black] node at (1.125, 3.375) {$0$};
		\draw [thin, black] node at (1.875, 3.375) {$1$};
		\draw [thin, black] node at (2.625, 3.375) {$0$};
		\draw [thin, black] node at (3.375, 3.375) {$0$};
		\draw [thin, black] node at (4.125, 3.375) {$0$};
		\draw [thin, black] node at (4.875, 3.375) {$\cdots$};

		\draw [thin, black] node at (1.875, .375) {$\vdots$};
		\draw [thin, black] node at (1.875, 1.125) {$0$};
		\draw [thin, black] node at (1.875, 1.875) {$0$};
		\draw [thin, black] node at (1.875, 2.625) {$0$};
		\draw [thin, black] node at (1.875, 4.125) {$0$};
		\draw [thin, black] node at (1.875, 4.875) {$\vdots$};

		\draw [black, very thick] (1.5,3.75)--(1.5,2.25)--(3,2.25)--(3,3.75)--cycle;
		\end{tikzpicture}
	}
	\subfigure{
		\begin{tikzpicture}
		\draw [step=.75, black, very thin] (0,0) grid (5.25,5.25);

		\draw [thin, black] node at (.375,3.375) {$\cdots$};
		\draw [thin, black] node at (1.125, 3.375) {$0$};
		\draw [thin, black] node at (1.875, 3.375) {\scriptsize$1\cdot1$};
		\draw [thin, black] node at (2.625, 3.375) {\scriptsize$0\cdot3$};
		\draw [thin, black] node at (3.375, 3.375) {$0$};
		\draw [thin, black] node at (4.125, 3.375) {$0$};
		\draw [thin, black] node at (4.875, 3.375) {$\cdots$};

		\draw [thin, black] node at (1.875, .375) {$\vdots$};
		\draw [thin, black] node at (1.875, 1.125) {$0$};
		\draw [thin, black] node at (1.875, 1.875) {$0$};
		\draw [thin, black] node at (1.875, 2.625) {\scriptsize$0\cdot2$};
		\draw [thin, black] node at (1.875, 4.125) {$0$};
		\draw [thin, black] node at (1.875, 4.875) {$\vdots$};

		\draw [thin, black] node at (2.625,2.625) {\scriptsize$z\cdot-1$};

		\draw [black, very thick] (1.5,3.75)--(1.5,2.25)--(3,2.25)--(3,3.75)--cycle;
		\end{tikzpicture}
	}
	\subfigure{
		\begin{tikzpicture}
		\draw [step=.75, black, very thin] (0,0) grid (5.25,5.25);

		\draw [thin, black] node at (.375,3.375) {$\cdots$};
		\draw [thin, black] node at (1.125, 3.375) {$0$};
		\draw [thin, black] node at (1.875, 3.375) {$1$};
		\draw [thin, black] node at (2.625, 3.375) {$0$};
		\draw [thin, black] node at (3.375, 3.375) {$0$};
		\draw [thin, black] node at (4.125, 3.375) {$0$};
		\draw [thin, black] node at (4.875, 3.375) {$\cdots$};

		\draw [thin, black] node at (1.875, .375) {$\vdots$};
		\draw [thin, black] node at (1.875, 1.125) {$0$};
		\draw [thin, black] node at (1.875, 1.875) {$0$};
		\draw [thin, black] node at (1.875, 2.625) {$0$};
		\draw [thin, black] node at (1.875, 4.125) {$0$};
		\draw [thin, black] node at (1.875, 4.875) {$\vdots$};

		\draw [thin, black] node at (2.625,2.625) {$1$};

		\draw [black, very thick] (1.5,3.75)--(1.5,2.25)--(3,2.25)--(3,3.75)--cycle;
		\end{tikzpicture}
	}
\end{figure}

%
%

Notice that if we look at the original recurrence $A_{r-1,c-1}+3A_{r-1,c}+2A_{r,c-1}-A_{r,c}=0$ with $r=1,c=1$, we have the same computation to solve for $A_{1,1}$ as when we use the overlay to compute $A_{1,1}$.  By shifting the overlay one position to the right, then $r=1, c=2$, which corresponds to the equation $0*1+0*3+1*2-z*-1=0$.  Solving for $z$, we have that $z=2$.

\begin{figure}[H]
\centering
	\subfigure{
		\begin{tikzpicture}
		\draw [step=.75, black, very thin] (0,0) grid (5.25,5.25);

		\draw [thin, black] node at (.375,3.375) {$\cdots$};
		\draw [thin, black] node at (1.125, 3.375) {$0$};
		\draw [thin, black] node at (1.875, 3.375) {$1$};
		\draw [thin, black] node at (2.625, 3.375) {$0$};
		\draw [thin, black] node at (3.375, 3.375) {$0$};
		\draw [thin, black] node at (4.125, 3.375) {$0$};
		\draw [thin, black] node at (4.875, 3.375) {$\cdots$};

		\draw [thin, black] node at (1.875, .375) {$\vdots$};
		\draw [thin, black] node at (1.875, 1.125) {$0$};
		\draw [thin, black] node at (1.875, 1.875) {$0$};
		\draw [thin, black] node at (1.875, 2.625) {$0$};
		\draw [thin, black] node at (1.875, 4.125) {$0$};
		\draw [thin, black] node at (1.875, 4.875) {$\vdots$};
		\draw [thin, black] node at (2.625,2.625) {$1$};

		\draw [black, very thick] (2.25,3.75)--(2.25,2.25)--(3.75,2.25)--(3.75,3.75)--cycle;
		\end{tikzpicture}
	}
	\subfigure{
		\begin{tikzpicture}
		\draw [step=.75, black, very thin] (0,0) grid (5.25,5.25);

		\draw [thin, black] node at (.375,3.375) {$\cdots$};
		\draw [thin, black] node at (1.125, 3.375) {$0$};
		\draw [thin, black] node at (1.875, 3.375) {$1$};
		\draw [thin, black] node at (2.625, 3.375) {\scriptsize$0\cdot1$};
		\draw [thin, black] node at (3.375, 3.375) {\scriptsize$0\cdot3$};
		\draw [thin, black] node at (4.125, 3.375) {$0$};
		\draw [thin, black] node at (4.875, 3.375) {$\cdots$};

		\draw [thin, black] node at (1.875, .375) {$\vdots$};
		\draw [thin, black] node at (1.875, 1.125) {$0$};
		\draw [thin, black] node at (1.875, 1.875) {$0$};
		\draw [thin, black] node at (1.875, 2.625) {$0$};
		\draw [thin, black] node at (2.625, 2.625) {\scriptsize$1\cdot2$};
		\draw [thin, black] node at (1.875, 4.125) {$0$};
		\draw [thin, black] node at (1.875, 4.875) {$\vdots$};

		\draw [thin, black] node at (3.375,2.625) {\scriptsize$z\cdot-1$};

		\draw [black, very thick] (2.25,3.75)--(2.25,2.25)--(3.75,2.25)--(3.75,3.75)--cycle;
		\end{tikzpicture}
	}
	\subfigure{
		\begin{tikzpicture}
		\draw [step=.75, black, very thin] (0,0) grid (5.25,5.25);

		\draw [thin, black] node at (.375,3.375) {$\cdots$};
		\draw [thin, black] node at (1.125, 3.375) {$0$};
		\draw [thin, black] node at (1.875, 3.375) {$1$};
		\draw [thin, black] node at (2.625, 3.375) {$0$};
		\draw [thin, black] node at (3.375, 3.375) {$0$};
		\draw [thin, black] node at (4.125, 3.375) {$0$};
		\draw [thin, black] node at (4.875, 3.375) {$\cdots$};

		\draw [thin, black] node at (1.875, .375) {$\vdots$};
		\draw [thin, black] node at (1.875, 1.125) {$0$};
		\draw [thin, black] node at (1.875, 1.875) {$0$};
		\draw [thin, black] node at (1.875, 2.625) {$0$};
		\draw [thin, black] node at (2.625, 2.625) {$1$};
		\draw [thin, black] node at (1.875, 4.125) {$0$};
		\draw [thin, black] node at (1.875, 4.875) {$\vdots$};

		\draw [thin, black] node at (3.375,2.625) {$2$};

		\draw [black, very thick] (2.25,3.75)--(2.25,2.25)--(3.75,2.25)--(3.75,3.75)--cycle;
		\end{tikzpicture}
	}
\end{figure}

By shifting the array one unit down from the first position of the overlay, this is another step that we could have taken as our second movement as it does not rely on the computation done above (although it will be left in the array).  For this computation, we solve for $z$ in the equation $0*1+0*3+1*2-z*-1=0$, which results in $z=3$.

\begin{figure}[H]
\centering
	\subfigure{
		\begin{tikzpicture}
		\draw [step=.75, black, very thin] (0,0) grid (5.25,5.25);

		\draw [thin, black] node at (.375,3.375) {$\cdots$};
		\draw [thin, black] node at (1.125, 3.375) {$0$};
		\draw [thin, black] node at (1.875, 3.375) {$1$};
		\draw [thin, black] node at (2.625, 3.375) {$0$};
		\draw [thin, black] node at (3.375, 3.375) {$0$};
		\draw [thin, black] node at (4.125, 3.375) {$0$};
		\draw [thin, black] node at (4.875, 3.375) {$\cdots$};

		\draw [thin, black] node at (1.875, .375) {$\vdots$};
		\draw [thin, black] node at (1.875, 1.125) {$0$};
		\draw [thin, black] node at (1.875, 1.875) {$0$};
		\draw [thin, black] node at (1.875, 2.625) {$0$};
		\draw [thin, black] node at (1.875, 4.125) {$0$};
		\draw [thin, black] node at (1.875, 4.875) {$\vdots$};
		\draw [thin, black] node at (2.625,2.625) {$1$};
		\draw [thin, black] node at (3.375,2.625) {$2$};

		\draw [black, very thick] (1.5,3)--(1.5,1.5)--(3,1.5)--(3,3)--cycle;
		\end{tikzpicture}
	}
	\subfigure{
		\begin{tikzpicture}
		\draw [step=.75, black, very thin] (0,0) grid (5.25,5.25);

		\draw [thin, black] node at (.375,3.375) {$\cdots$};
		\draw [thin, black] node at (1.125, 3.375) {$0$};
		\draw [thin, black] node at (1.875, 3.375) {$1$};
		\draw [thin, black] node at (2.625, 3.375) {$0$};
		\draw [thin, black] node at (3.375, 3.375) {$0$};
		\draw [thin, black] node at (4.125, 3.375) {$0$};
		\draw [thin, black] node at (4.875, 3.375) {$\cdots$};

		\draw [thin, black] node at (1.875, .375) {$\vdots$};
		\draw [thin, black] node at (1.875, 1.125) {$0$};
		\draw [thin, black] node at (1.875, 1.875) {\scriptsize$0\cdot2$};
		\draw [thin, black] node at (1.875, 2.625) {\scriptsize$0\cdot1$};
		\draw [thin, black] node at (2.625, 2.625) {\scriptsize$1\cdot3$};
		\draw [thin, black] node at (1.875, 4.125) {$0$};
		\draw [thin, black] node at (1.875, 4.875) {$\vdots$};
		\draw [thin, black] node at (3.375,2.625) {$2$};

		\draw [thin, black] node at (2.625,1.875) {\scriptsize$z\cdot-1$};

		\draw [black, very thick] (1.5,3)--(1.5,1.5)--(3,1.5)--(3,3)--cycle;
		\end{tikzpicture}
	}
	\subfigure{
		\begin{tikzpicture}
		\draw [step=.75, black, very thin] (0,0) grid (5.25,5.25);

		\draw [thin, black] node at (.375,3.375) {$\cdots$};
		\draw [thin, black] node at (1.125, 3.375) {$0$};
		\draw [thin, black] node at (1.875, 3.375) {$1$};
		\draw [thin, black] node at (2.625, 3.375) {$0$};
		\draw [thin, black] node at (3.375, 3.375) {$0$};
		\draw [thin, black] node at (4.125, 3.375) {$0$};
		\draw [thin, black] node at (4.875, 3.375) {$\cdots$};

		\draw [thin, black] node at (1.875, .375) {$\vdots$};
		\draw [thin, black] node at (1.875, 1.125) {$0$};
		\draw [thin, black] node at (1.875, 1.875) {$0$};
		\draw [thin, black] node at (1.875, 2.625) {$0$};
		\draw [thin, black] node at (2.625, 2.625) {$1$};
		\draw [thin, black] node at (1.875, 4.125) {$0$};
		\draw [thin, black] node at (1.875, 4.875) {$\vdots$};
		\draw [thin, black] node at (3.375,2.625) {$2$};

		\draw [thin, black] node at (2.625,1.875) {$3$};

		\draw [black, very thick] (1.5,3)--(1.5,1.5)--(3,1.5)--(3,3)--cycle;
		\end{tikzpicture}
	}
\end{figure}

By shifting the overlay in different positions throughout $A$ such that three of the four corresponding entries in $A$ are known, we can populate the rest of $A$, as seen below.

\begin{center}
\begin{tikzpicture}
		\draw [step=.75, black, very thin] (0,0) grid (5.25,5.25);

		\draw [thin, black] node at (.375, .45) {\reflectbox{$\ddots$}};
		\draw [thin, black] node at (1.125, .45) {$\vdots$};
		\draw [thin, black] node at (1.875, .45) {$\vdots$};
		\draw [thin, black] node at (2.625, .45) {$\vdots$};
		\draw [thin, black] node at (3.375, .45) {$\vdots$};
		\draw [thin, black] node at (4.125, .45) {$\vdots$};
		\draw [thin, black] node at (4.875, .45) {$\ddots$};

		\draw [thin, black] node at (.375, 1.125) {$\cdots$};
		\draw [thin, black] node at (1.125, 1.125) {$-\frac{3}{8}$};
		\draw [thin, black] node at (1.875, 1.125) {$0$};
		\draw [thin, black] node at (2.625, 1.125) {$9$};
		\draw [thin, black] node at (3.375, 1.125) {$58$};
		\draw [thin, black] node at (4.125, 1.125) {$249$};
		\draw [thin, black] node at (4.875, 1.125) {$\cdots$};

		\draw [thin, black] node at (.375, 1.875) {$\cdots$};
		\draw [thin, black] node at (1.125, 1.875) {$\frac{3}{4}$};
		\draw [thin, black] node at (1.875, 1.875) {$0$};
		\draw [thin, black] node at (2.625, 1.875) {$3$};
		\draw [thin, black] node at (3.375, 1.875) {$13$};
		\draw [thin, black] node at (4.125, 1.875) {$40$};
		\draw [thin, black] node at (4.875, 1.875) {$\cdots$};

		\draw [thin, black] node at (.375, 2.625) {$\cdots$};
		\draw [thin, black] node at (1.125, 2.625) {$-\frac{3}{2}$};
		\draw [thin, black] node at (1.875, 2.625) {$0$};
		\draw [thin, black] node at (2.625, 2.625) {$1$};
		\draw [thin, black] node at (3.375, 2.625) {$2$};
		\draw [thin, black] node at (4.125, 2.625) {$4$};
		\draw [thin, black] node at (4.875, 2.625) {$\cdots$};

		\draw [thin, black] node at (.375,3.375) {$\cdots$};
		\draw [thin, black] node at (1.125, 3.375) {$0$};
		\draw [thin, black] node at (1.875, 3.375) {$1$};
		\draw [thin, black] node at (2.625, 3.375) {$0$};
		\draw [thin, black] node at (3.375, 3.375) {$0$};
		\draw [thin, black] node at (4.125, 3.375) {$0$};
		\draw [thin, black] node at (4.875, 3.375) {$\cdots$};

		\draw [thin, black] node at (.375, 4.125) {$\cdots$};
		\draw [thin, black] node at (1.125, 4.125) {$1$};
		\draw [thin, black] node at (1.875, 4.125) {$0$};
		\draw [thin, black] node at (2.625, 4.125) {$-\frac{2}{3}$};
		\draw [thin, black] node at (3.375, 4.125) {$\frac{2}{9}$};
		\draw [thin, black] node at (4.125, 4.125) {$-\frac{2}{27}$};
		\draw [thin, black] node at (4.875, 4.125) {$\cdots$};

		\draw [thin, black] node at (.375, 4.95) {$\ddots$};
		\draw [thin, black] node at (1.125, 4.95) {$\vdots$};
		\draw [thin, black] node at (1.875, 4.95) {$\vdots$};
		\draw [thin, black] node at (2.625, 4.95) {$\vdots$};
		\draw [thin, black] node at (3.375, 4.95) {$\vdots$};
		\draw [thin, black] node at (4.125, 4.95) {$\vdots$};
		\draw [thin, black] node at (4.875, 4.95) {\reflectbox{$\ddots$}};

\end{tikzpicture}
\end{center}

%
%

This example illustrates the use of the overlay to populate $A$ rather than working directly with the recurrence. \\


\begin{theorem} 
\label{std_init_conditions}
If we have an overlay $B=\left(b_{i,j}\right)$ such that the top row of the overlay is $u+1$ units wide, the bottom row is $l+1$ units wide, the furthest column to the right containing a nonzero entry in the overlay is indexed zero, the furthest column to the left containing a nonzero entry in the overlay is indexed $n$, and the bottom-most row containing a nonzero entry in the overlay is indexed $0$ while the top-most row containing a nonzero entry in the overlay is indexed $m$, then $m$ consecutive rows of initial conditions, along with $u$ consecutive columns coming out of the top and $l$ consecutive columns coming out of the bottom of these rows is sufficient to populate the rest of the array.
\end{theorem}

\proof

Observe that if $n=0$ and $m=0$ then $b_{0,0}  \ne 0$ (otherwise we would not have a valid overlay).  Then the overlay corresponds to the recurrence
$$0=\sum_{i=0}^0 \sum_{j=0}^0 b_{i,j} A_{r-i,c-j}=b_{0,0} A_{r,c}.$$

Since $b_{0,0}\ne 0$, then $A_{r,c}=0$, so the only valid array is the zero array.  Suppose that $n$ and $m$ are not both simultaneously zero.

If $m=0$ and $n\ne 0$, then the initial conditions we claim to be sufficient to populate the entire array is $n$ columns that are pre-populated and no rows that are pre-populated.  Without loss of generality, let the pre-populated columns be indexed 0 through $n-1$.  Then we have $$\sum_{j=0}^n b_{0,j} A_{r,c-j}=0.$$

If trying to populate a column one unit to the right of the previously constructed columns, say column $c\ge n$, then by our assumptions, the columns indexed 0 through $c-1$ have already been populated.  Thus for our summation, we have
$$\sum_{j=0}^n b_{0,j} A_{r,c-j}=0.$$  Since $j$ ranges between 0 and $n$, then $c-1-j \le c-1$ and $0 \le c-n \le c-j$.  Thus each $A_{r,c-j}$ is already determined except for $A_{r,c}$, and since $b_{0,0} \ne 0$, $A_{r,c}$ is uniquely determined by the recurrence.  Since $r$ is arbitrary, then the entire column can be populated this way. 

Similarly, relying on the fact that $b_{0,n}\ne 0$, we can constrcut a column $c<0$ 

For the remaining cases, it suffices to show that we can construct a row above the $m$ consecutive rows and a row below the $m$ consecutive rows in order to populate the array because if we do one, then we can build the next row above by placing the overlay in the starting position from the previous row, move the overlay up one row, and repeat the process outlined below.  The same holds for rows below, we just move the overlay down instead of up.

Given the overlay described above, the overlay corresponds to the recurrence
$$ \sum_{i=0}^m \sum_{j=0}^n b_{i,j} A_{r-i,c-j}=0$$
where $b_{i,j}$ are the coefficients in the overlay with at least one nonzero entry in the $0^{\text{th}}$ row, in the $0^{\text{th}}$ column, in the $m^{\text{th}}$ row and in the $n^{\text{th}}$ column.

Without loss of generality, let the $m$ consecutive rows of initial conditions be the rows indexed 0 through $m-1$.

Case 1 ($u=0$): Since $u=0$, then there is a single nonzero value in the top row of the overlay. Let the nonzero element in the $m^{\text{th}}$ row of the overlay (top row) be the column indexed $a\in \Z$. Through the description of the overlay, we know that $0 \le s\le n$. Since there is a single nonzero value in the top row of the overlay, there are no columns of initial conditions protruding out of the top of the $m$ consecutive rows. If we place the bottom row of the overlay in the $m-1$ row, then the top row of the overlay is in the row indexed $-1$. Let $y\in \Z$. Since there are no columns of initial conditions, let the column in which the $s$ column in the overlay resides be indexed $y$. By the recurrence, we have
$$ 0=\sum_{i=0}^m \sum_{j=0}^n b_{i,j} A_{r-i,c-j}. $$


Taking $c=s+y$ and $r=m-1$ as described in the placement of the bottom row, our recurrence becomes
\begin{align*}
0 &= \sum_{i=0}^m \sum_{j=0}^n b_{i,j} A_{m-1-i, s+y-j} \\
0 &= b_{m,s} A_{m-1-m, s+y-s} + \sum_{i=0}^{m-1} \sum_{j=0}^{s-1} b_{i,j} A_{m-1-i, s+y-j} + \sum_{i=0}^{m-1} \sum_{j=s+1}^n b_{i,j} A_{m-1-i,s+y-j}\\
& \hspace{1 cm} +\sum_{i=0}^{m-1} b_{i,s} A_{m-1-i,s+y-s} + \sum_{j=0}^{s-1} b_{m,j} A_{m-1-m, s+y-j} + \sum_{j=s+1}^n b_{m,j} A_{m-1-m, s+y-j}
\end{align*}

Since the $m^{\text{th}}$ row of the overlay contains only one nonzero element, then $b_{m,j}=0$ for all $j\ne s$. Thus
$$\sum_{j=0}^{s-1} b_{m,j} A_{m-1-m, s+y-j} + \sum_{j=s+1}^{n} b_{m,j} A_{m-1-m, s+y-j}=0.$$

So simplifying some of the indices in the previous expansion, we have
\begin{align*}
0 &= b_{m,s} A_{-1, y} + \sum_{i=0}^{m-1} \sum_{j=0}^{s-1} b_{i,j} A_{m-1-i, s+y-j} + \sum_{i=0}^{m-1} \sum_{j=s+1}^n b_{i,j} A_{m-1-i,s+y-j}+\sum_{i=0}^{m-1} b_{i,s} A_{m-1-i,y} \\
&\Leftrightarrow -b_{m,s} A_{-1,y}= \sum_{i=0}^{m-1} \sum_{j=0}^{s-1} b_{i,j} A_{m-1-i, s+y-j} + \sum_{i=0}^{m-1} \sum_{j=s+1}^n b_{i,j} A_{m-1-i,s+y-j}+\sum_{i=0}^{m-1} b_{i,s} A_{m-1-i,y}.
\end{align*}

Since $b_{m,s} \ne 0$, dividing both sides by $-b_{m,s}$, we have
$$ A_{-1,y}= -\dfrac{1}{b_{m,s}} \left( \sum_{i=0}^{m-1} \sum_{j=0}^{s-1} b_{i,j} A_{m-1-i, s+y-j} + \sum_{i=0}^{m-1} \sum_{j=s+1}^n b_{i,j} A_{m-1-i,s+y-j}+\sum_{i=0}^{m-1} b_{i,s} A_{m-1-i,y} \right).$$

Since the rows in each sum range from $m-1-(m-1)=0$ to $m-1-0=m-1$, then each $A_{m-1-i,s+y-j}$ is determined by the initial conditions. Thus, $A_{-1,y}$ is uniquely determined by the recurrence. Since $y$ is arbitrary, then the entire row can be populated this way.\\

Case 2 ($u>0$): Let $a,d,\in \Z$. Then let the $u$ consecutive columns above the $m$ consecutive rows be the columns indexed $a$ through $a+u-1$ and the $l$ columns below the $m$ consecutive rows be indexed $d$ through $d+l-1$. Below, an array with the initial conditions present can be seen.

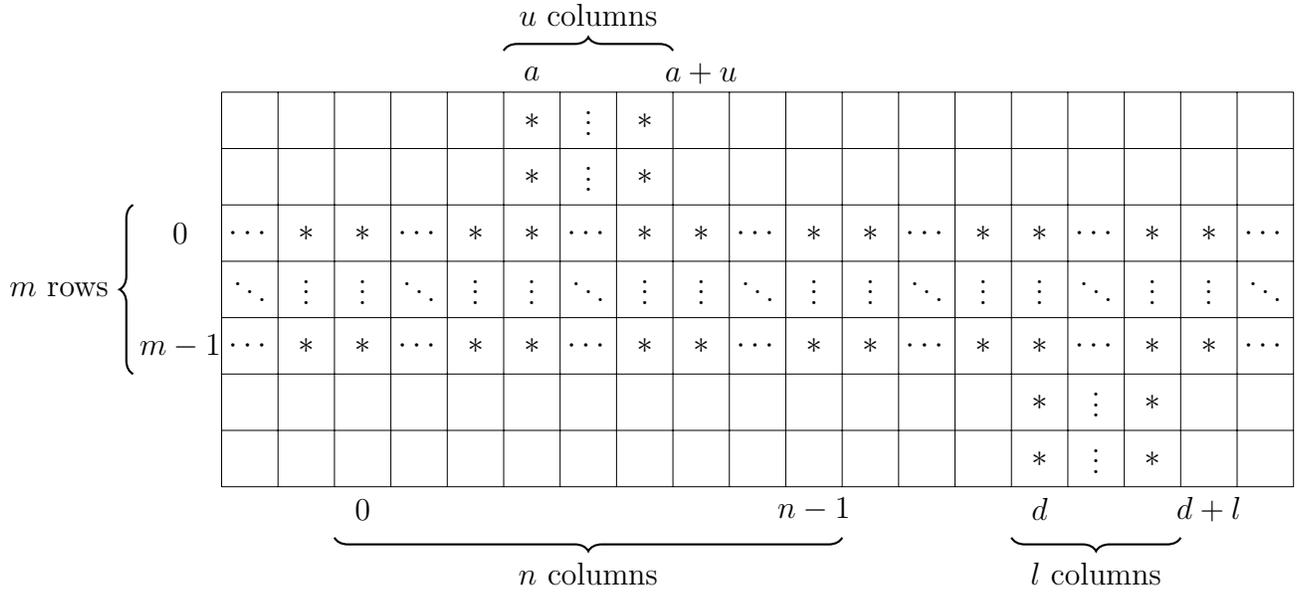
\begin{center}
\begin{tikzpicture}
\draw [step=.75, black, very thin] (0,0) grid (14.25,5.25);


\foreach \i in {.75,1.5,3,3.75,5.25,6,7.5,8.25,9.75,10.5,12,12.75} {
	\foreach \j in {1.5,3} {
		\draw [thin,black] node at (\i+.375, \j+.375) {$*$};
	}
	\draw [thin,black] node at (\i+.375,2.7) {$\vdots$};
}

\foreach \i in {0,2.25,4.5,6.75,9,11.25,13.5} {
	\foreach \j in {1.5,3} {
		\draw [thin,black] node at (\i+.375, \j+.375) {$\cdots$};
	}
	\draw [thin, black] node at (\i+.375, 2.7) {$\ddots$};
}

\foreach \i in {3.75, 5.25} {
	\foreach \j in {3.75, 4.5} {
		\draw [thin, black] node at (\i+.375, \j+.375) {$*$};
	}
}
\draw [thin, black] node at (4.875, 4.2) {$\vdots$};
\draw [thin, black] node at (4.875, 4.95) {$\vdots$};

\foreach \i in {10.5, 12} {
	\foreach \j in {0, .75} {
		\draw [thin, black] node at (\i+.375, \j+.375) {$*$};
	}
}
\draw [thin, black] node at (11.625, 0.45) {$\vdots$};
\draw [thin, black] node at (11.625, 1.2) {$\vdots$};

\draw [addbrace] (-1,1.5) -- node[left=10pt]{$m$ rows}(-1,3.75);
\draw [thin, black] node at (-.55, 3.375) {0};
\draw [thin, black] node at (-.55, 1.875){$m-1$};

\draw [addbrace] (3.75,5.625) -- node[above=10pt]{$u$ columns}(6,5.625);
\draw [thin, black] node at (4.125, 5.5) {$a$};
\draw [thin, black] node at (6.375, 5.5) {$a+u$};

\draw [addbrace] (12.75,-.5)--node[below=10pt]{$l$ columns}(10.5,-.5);
\draw [thin, black] node at (10.875, -.3) {$d$};
\draw [thin, black] node at (13.125, -.3) {$d+l$};

\draw [addbrace] (8.25,-.5)--node[below=10pt]{$n$ columns}(1.5,-.5);
\draw [thin, black] node at (7.875, -.3) {$n-1$};
\draw [thin, black] node at (1.875, -.3) {0};


\end{tikzpicture}
\end{center}

Observe that $u\le n$ and $l\le n$. Let the column of the right-most nonzero entry in the top row ($m^{\text{th}}$ row) in the overlay be denoted $s$ and let the column of the right-most nonzero entry in the bottom row ($0^{\text{th}}$ row) in the overlay be denoted $t$. If we place the bottom row of the overlay in the $m-1$ row of the array, then the top row is in the row indexed $-1$. If we then slide the overlay horizontally so that the left-most nonzero entry in the top row of the overlay is in the column indexed by $a$ in the array. By positioning the overlay in this manner, this corresponds to taking $c=s+a+u$ and $r=m-1$ in the recurrence, as described in the placement of the bottom row, as demonstrated in the image below.

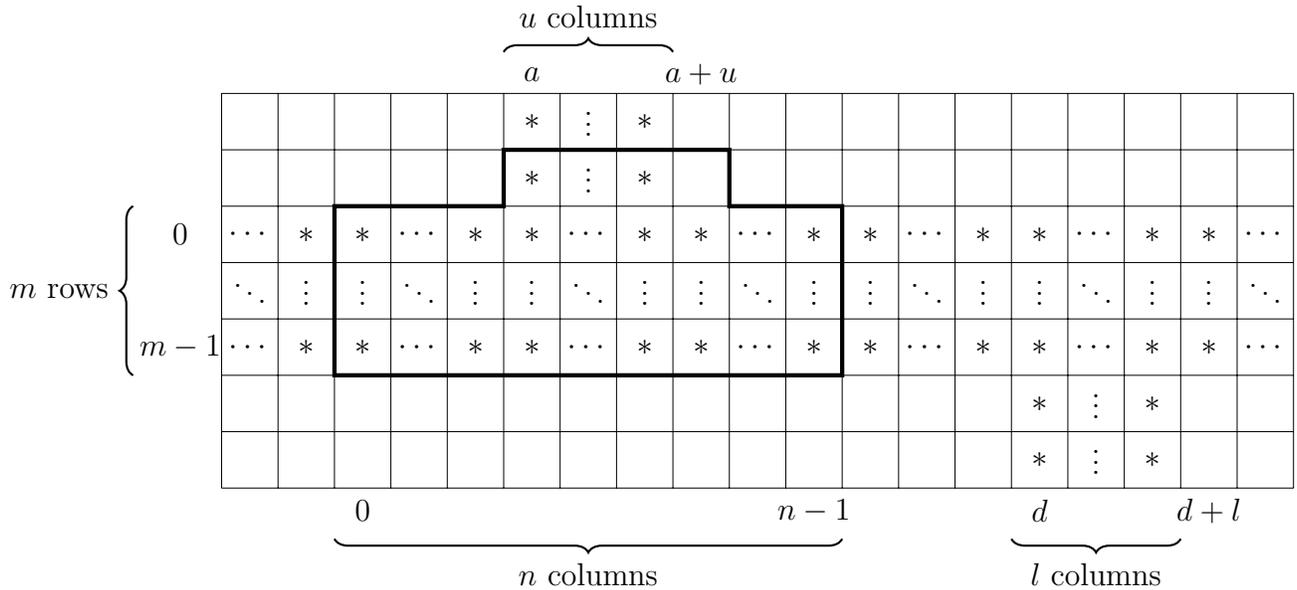
\begin{center}
\begin{tikzpicture}
\draw [step=.75, black, very thin] (0,0) grid (14.25,5.25);


\foreach \i in {.75,1.5,3,3.75,5.25,6,7.5,8.25,9.75,10.5,12,12.75} {
	\foreach \j in {1.5,3} {
		\draw [thin,black] node at (\i+.375, \j+.375) {$*$};
	}
	\draw [thin,black] node at (\i+.375,2.7) {$\vdots$};
}

\foreach \i in {0,2.25,4.5,6.75,9,11.25,13.5} {
	\foreach \j in {1.5,3} {
		\draw [thin,black] node at (\i+.375, \j+.375) {$\cdots$};
	}
	\draw [thin, black] node at (\i+.375, 2.7) {$\ddots$};
}

\foreach \i in {3.75, 5.25} {
	\foreach \j in {3.75, 4.5} {
		\draw [thin, black] node at (\i+.375, \j+.375) {$*$};
	}
}
\draw [thin, black] node at (4.875, 4.2) {$\vdots$};
\draw [thin, black] node at (4.875, 4.95) {$\vdots$};

\foreach \i in {10.5, 12} {
	\foreach \j in {0, .75} {
		\draw [thin, black] node at (\i+.375, \j+.375) {$*$};
	}
}
\draw [thin, black] node at (11.625, 0.45) {$\vdots$};
\draw [thin, black] node at (11.625, 1.2) {$\vdots$};

\draw [black, ultra thick] (3.75,4.5)--(6.75,4.5)--(6.75,3.75)--(8.25,3.75)--(8.25,1.5)--(1.5,1.5)--(1.5,3.75)--(3.75,3.75)--cycle;

\draw [addbrace] (-1,1.5) -- node[left=10pt]{$m$ rows}(-1,3.75);
\draw [thin, black] node at (-.55, 3.375) {0};
\draw [thin, black] node at (-.55, 1.875){$m-1$};

\draw [addbrace] (3.75,5.625) -- node[above=10pt]{$u$ columns}(6,5.625);
\draw [thin, black] node at (4.125, 5.5) {$a$};
\draw [thin, black] node at (6.375, 5.5) {$a+u$};

\draw [addbrace] (12.75,-.5)--node[below=10pt]{$l$ columns}(10.5,-.5);
\draw [thin, black] node at (10.875, -.3) {$d$};
\draw [thin, black] node at (13.125, -.3) {$d+l$};

\draw [addbrace] (8.25,-.5)--node[below=10pt]{$n$ columns}(1.5,-.5);
\draw [thin, black] node at (7.875, -.3) {$n-1$};
\draw [thin, black] node at (1.875, -.3) {0};


\end{tikzpicture}
\end{center}

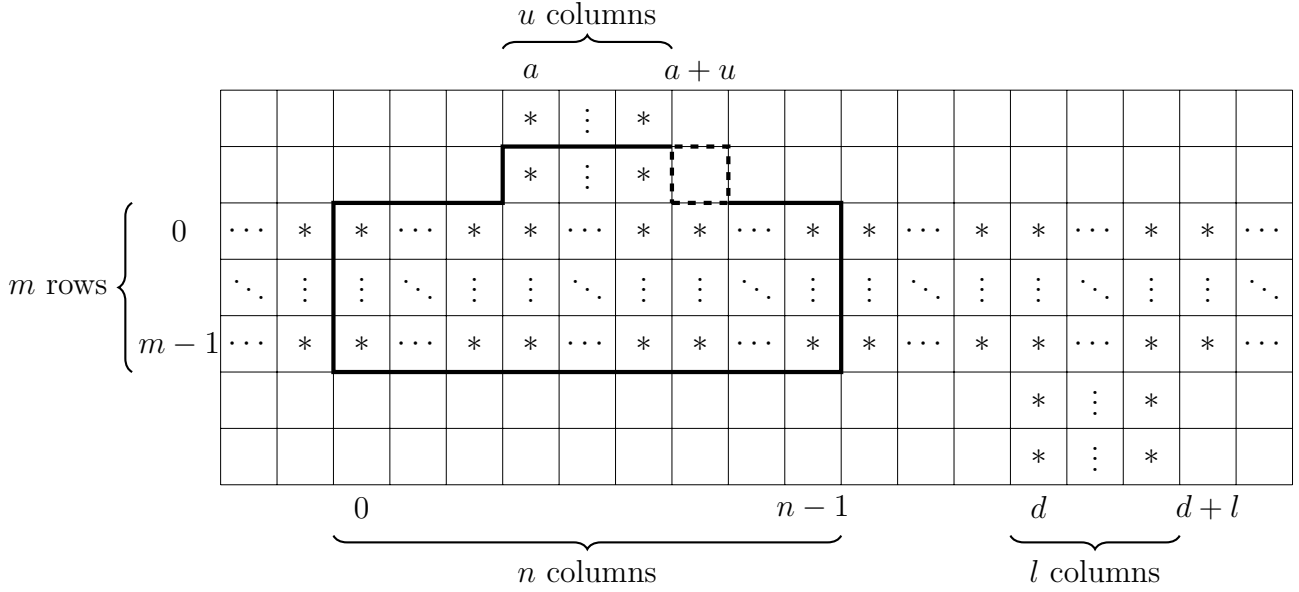
\begin{center}
\begin{tikzpicture}
\draw [step=.75, black, very thin] (0,0) grid (14.25,5.25);


\foreach \i in {.75,1.5,3,3.75,5.25,6,7.5,8.25,9.75,10.5,12,12.75} {
	\foreach \j in {1.5,3} {
		\draw [thin,black] node at (\i+.375, \j+.375) {$*$};
	}
	\draw [thin,black] node at (\i+.375,2.7) {$\vdots$};
}

\foreach \i in {0,2.25,4.5,6.75,9,11.25,13.5} {
	\foreach \j in {1.5,3} {
		\draw [thin,black] node at (\i+.375, \j+.375) {$\cdots$};
	}
	\draw [thin, black] node at (\i+.375, 2.7) {$\ddots$};
}

\foreach \i in {3.75, 5.25} {
	\foreach \j in {3.75, 4.5} {
		\draw [thin, black] node at (\i+.375, \j+.375) {$*$};
	}
}
\draw [thin, black] node at (4.875, 4.2) {$\vdots$};
\draw [thin, black] node at (4.875, 4.95) {$\vdots$};

\foreach \i in {10.5, 12} {
	\foreach \j in {0, .75} {
		\draw [thin, black] node at (\i+.375, \j+.375) {$*$};
	}
}
\draw [thin, black] node at (11.625, 0.45) {$\vdots$};
\draw [thin, black] node at (11.625, 1.2) {$\vdots$};

\draw [black, ultra thick] (6.75,3.75)--(8.25,3.75)--(8.25,1.5)--(1.5,1.5)--(1.5,3.75)--(3.75,3.75)--(3.75,4.5)--(6,4.5);
\draw [black, ultra thick, dashed] (6,4.5)--(6.75,4.5)--(6.75,3.75)--(6,3.75)--cycle;

\draw [addbrace] (-1,1.5) -- node[left=10pt]{$m$ rows}(-1,3.75);
\draw [thin, black] node at (-.55, 3.375) {0};
\draw [thin, black] node at (-.55, 1.875){$m-1$};

\draw [addbrace] (3.75,5.625) -- node[above=10pt]{$u$ columns}(6,5.625);
\draw [thin, black] node at (4.125, 5.5) {$a$};
\draw [thin, black] node at (6.375, 5.5) {$a+u$};

\draw [addbrace] (12.75,-.5)--node[below=10pt]{$l$ columns}(10.5,-.5);
\draw [thin, black] node at (10.875, -.3) {$d$};
\draw [thin, black] node at (13.125, -.3) {$d+l$};

\draw [addbrace] (8.25,-.5)--node[below=10pt]{$n$ columns}(1.5,-.5);
\draw [thin, black] node at (7.875, -.3) {$n-1$};
\draw [thin, black] node at (1.875, -.3) {0};


\end{tikzpicture}
\end{center}

Thus our recurrence becomes

\begin{align*}
0 &= \sum_{i=0}^m \sum_{j=0}^n b_{i,j} A_{m-1-i, s+a+u-j}\\
\Leftrightarrow 0 &= b_{m,s} A_{m-1-m, s+a-u-s} + \sum_{i=0}^{m-1} \sum_{j=0}^{s-1} b_{i,j} A_{m-1-i, s+a+u-j} + \sum_{i=0}^{m-1} \sum_{j=s+1}^n b_{i,j} A_{m-1-i,s+a+u-j}\\
& \hspace{1 cm} +\sum_{i=0}^{m-1} b_{i,s} A_{m-1-i,s+a+u-s} + \sum_{j=0}^{s-1} b_{m,j} A_{m-1-m, s+a+u-j} + \sum_{j=s+1}^n b_{m,j} A_{m-1-m, s+a+u-j}
\end{align*}

Since the $m^{\text{th}}$ row of the overlay contains only entries from $s$ to $s+u$, then $b_{m,j}=0$ for all $j<s$, and $b_{m,j}=0$ for all $j>s+u$. Thus,
$$\sum_{j=0}^{s-1} b_{m,j}A_{m-1-m,s+a+u-j}=0+\sum_{j=s+1}^{s+u}b_{m,j}A_{m-1-m,s+a+u-j}.$$


So simplifying some of our indices from the previous expansion, we have

\begin{align*}
0=b_{m,s}A_{-1,a+u} &+ \sum_{i=0}^{m-1}\sum_{j=0}^{s-1}b_{i,j}A_{m-1-i,s+a+u-j} + \sum_{i=0}^{m-1}\sum_{j=s+1}^{n}b_{i,j}A_{m-1-i,s+a+u-j}\\
&+\sum_{i=0}^{m-1}b_{i,s}A_{m-1-i,a+u} +\sum_{j=s+1}^{s+u}b_{m,j}A_{-1,s+a+u-j} \\
\Leftrightarrow -b_{m,s}A_{-1,a+u} &=\sum_{i=0}^{m-1} \sum_{j=0}^{s-1}b_{i,j}A_{m-1-i,s+a+u-j}+\sum_{i=0}^{m-1}\sum_{j=s+1}^n b_{i,j} A_{m-1-i,s+a+u-j}\\
&+\sum_{i=0}^{m-1}b_{i,s}A_{m-1-i,a+u} +\sum_{j=s+1}^{s+u}b_{m,j}A_{-1,s+a+u-j}.
\end{align*}

Since $b_{m,s} \ne 0$, dividing both sides by $-b_{m,s}$, we have

\begin{align*}
A_{-1,a+u}&=-\dfrac{1}{b_{m,s}} \left( \sum_{i=0}^{m-1} \sum_{j=0}^{s-1}b_{i,j}A_{m-1-i,s+a+u-j} +\sum_{i=0}^{m-1} \sum_{j=s+1}^n b_{i,j}A_{m-1-i,s+a+u-j} \right. \\
&+\left. \sum_{i=0}^{m-1}b_{i,s}A_{m-1-i,a+u} + \sum_{j=s+1}^{s+u}b_{m,j}A_{-1,s+a+u-j} \right) .
\end{align*}

Since the rows of the array in each of the first three sums range from $m-1-(m-1)=0$ to $m-1-0=m-1$, then each $A_{m-1-i,s+a+u-j}$ is determined by the initial conditions. Now let's consider the $m^{\text{th}}$ row of the overlay, which is denoted by the last sum and is the row indexed $-1$ in the array. In the fourth sum, the columns in the array range from $s+a+u-(s+u)=a$ to $s+a+u-(s+1)=a+u-1$. Since the columns in which the initial conditions are defined are columns $a$ through $a+u-1$, then each $A_{-1,s+a+u-j}$ is determined by the initial conditions.

Since all entries on the right side of the equals sign are determined by the initial conditions, then $A_{-1,a+u}$ is uniquely determined by the recurrence. By shifting the overlay to the right one column, we can similarly get the next entry and thus fill in the right side of the row. If we orient the overlay in the same location vertically but shift the overlay one column to the left compared to the original orientation, then this corresponds to $r=m-1$ and $c=s+a+u-1$ in the recurrence. Thus our recurrence becomes

\begin{align*}
0 &= \sum_{i=0}^m \sum_{j=0}^n b_{i,j} A_{m-1-i, s+u-1-j} \\
\Leftrightarrow 0 &= b_{m,s+u}A_{m-1-m, s+a+u-1-(s+u)} + \sum_{i=0}^{m-1} \sum_{j=0}^{s+u-1} b_{i,j} A_{m-1-i, s+a+u-1-j} \\
&\hspace{1 cm} + \sum_{i=0}^{m-1} \sum_{j=s+u+1}^n b_{i,j} A_{m-1-i, s+a+u-1-j} + \sum_{i=0}^{m-1} b_{i,s+u} A_{m-1-i, s+a+u-1-(s+u)} \\
&\hspace{1 cm} + \sum_{j=0}^{s+u-1} b_{m,j} A_{m-1-m, s+a+u-1-j} + \sum_{j=s+u+1}^n b_{m,j} A_{m-1-m, s+a+u-1-j}
\end{align*}

Since the $m^{\text{th}}$ row of the overlay contains only entries from $s$ to $s+u$, then $b_{m,j}=0$ for all $j<s$ and $b_{m,j}=0$ for all $j>s+u$. Thus
$$\sum_{j=0}^{s+u-1} b_{m,j}A_{m-1-m, s+a+u-1-j} + \sum_{j=s+u+1}^n b_{m,j} A_{m-1-m, s+a+u-1-j} = \sum_{j=s}^{s+u-1} b_{m,j} A_{m-1-m, s+a+u-1-j} + 0. $$


So simplifying some of the indices in the previous expansion, we have
\begin{align*}
0 &= b_{m,s+u} A_{-1,a-1} + \sum_{i=0}^{m-1} \sum_{j=0}^{s+u-1} b_{i,j} A_{m-1-i, s+a+u-1-j} + \sum_{i=0}^{m-1} \sum_{j=s+u+1}^n b_{i,j} A_{m-1-i, s+a+u-1-j} \\
& \hspace{1 cm} +\sum_{i=0}^{m-1} b_{i,s+u} A_{m-1-i, a-1} + \sum_{j=s}^{s+u-1} b_{m,j} A_{-1, s+a+u-1-j} \\
&\Leftrightarrow -b_{m,s+u}A_{-1,a-1} = \sum_{i=0}^{m-1} \sum_{j=0}^{s+u-1}b_{i,j} A_{m-1-i, s+a+u-1-j} + \sum_{i=0}^{m-1} \sum_{j=s+u+1}^n b_{i,j} A_{m-1-i, s+a+u-1-j} \\
&\hspace{1 cm} + \sum_{i=0}^{m-1} b_{i,s+u} A_{m-1-i,a-1} + \sum_{j=s}^{s+u-1} b_{m,j} A_{-1,s+a+u-1-j}.
\end{align*}

Since $b_{m,s+u} \ne 0$, dividing both sides by $-b_{m,s+u}$, we have
\begin{align*}
A_{-1,a-1} &= -\dfrac{1}{b_{m,s+u}} \left( \sum_{i=0}^{m-1} \sum_{j=0}^{s+u-1} b_{i,j} A_{m-1-i, s+a+u-1-j} + \sum_{i=0}^{m-1} \sum_{j=s+u+1}^n b_{i,j} A_{m-1-i, s+a+u-1-j} \right. \\
& \hspace{1 cm} + \left. \sum_{i=0}^{m-1} b_{i,s+u} A_{m-1-i, a-1} + \sum_{j=s}^{s+u-1}b_{m,j} A_{-1, s+a+u-1-j} \right) .
\end{align*}

Since the rows of the array in each of the first three sums range from $m-1-(m-1)=0$ to $m-1-0=m-1$, then each $A_{m-1-i, s+a+u-1-j}$ is determined by the initial conditions. Now let's consider the $m^\text{th}$ row of the overlay, which is denoted by the last sum and is the row indexed $-1$ in the array. In the fourth sum, the columns in the array range from $s+a+u-1-(s+u-1)=a$ to $s+a+u-1-(s)=a+u-1.$ Since the columns in which the initial conditions are defined are columns $a$ through $a+u-1$, then each $A_{-1,s+a+u-1-j}$ is determined by the initial conditions.

Since all entries on the right side of the equals sign are determined by the initial conditions, then $A_{-1,a-1}$ is uniquely determined by the recurrence.

By shifting the overlay horizontally to the left by one column, then we can similarly get the next entry, and thus fill in the left side of the first row. In order to get the row above this, we shift the overlay up one row from its initial position and repeat this process, thus filling in the negative rows of the array.

For the positive rows indexed $\ge m$, the proof is similar, which $l$ is used instead of $u$, $d$ is used instead of $a$, $t$ is used instead of $s$, $m$ is used instead of $-1$, and a shift of the overlay up a row becomes a shift of the overlay down a row.  Alternatively, we can think about getting the positive rows through plane transformations.

\newpage

\begin{definition}
A \emph{valid initial conditions layout} is the set of elements in $\Z^2$ (coordinates) that define the location of the initial conditions.
\end{definition}

In order to be an initial conditions layout, the initial conditions must construct a unique array when the overlay is applied.  In other words, a layout is the locations on the array that must be filled in with initial values in order to construct a unique array for a given overlay.  We do assume here that the overlay contains at least one nonzero entry.

\begin{definition}
A \emph{standard initial conditions layout} is the layout described in theorem \ref{std_init_conditions}.
\end{definition}

However, there are other possible initial conditions layouts that are valid for different overlays.  For example, if we have the template $T=b_{0,0}+b_{1,0}Y+b_{1,1}XY$, where $b_{0,0}, b_{1,0}, b_{1,1} \ne 0$, with corresponding overlay

\begin{center}
\begin{tikzpicture}
\draw [thin, black] (1.5,0)--(1.5,1.5)--(0,1.5)--(0,.75)--(.75,.75)--(.75,0)--cycle;
\draw [thin, black] (1.5,.75)--(.75,.75)--(.75,1.5);

\draw [thin, black] node at (.375, 1.125) {$b_{1,1}$};
\draw [thin, black] node at (1.125, 1.125) {$b_{1,0}$};
\draw [thin, black] node at (1.125, .375) {$b_{0,0}$};

\end{tikzpicture}
\end{center}

\noindent then the standard initial conditions layout would be

\begin{center}
\begin{tikzpicture}
\draw [step=.75, black, very thin] (0,0) grid (6.75,6.75);

\foreach \i in {.75, 1.5, 2.25, 3, 3.75, 4.5, 5.25} {
	\draw [thin, black] node at (\i+.375, 1.875) {$*$};
	}

\foreach \i in {2.25, 3, 3.75, 4.5, 5.25} {
	\draw [thin, black] node at (3.425, \i+.375) {$*$};
	}

\foreach \i in {0, 6} {
	\draw [thin, black] node at (\i+.375, 1.875) {$\cdots$};
	}

\draw [thin, black] node at (3.425, 6.375) {$\vdots$};

\end{tikzpicture}
\end{center}

\noindent where * denotes a populated cell due to an initial condition. Alternatively, another valid initial conditions layout could be

\begin{center}
\begin{tikzpicture}
\draw [step=.75, black, very thin] (0,0) grid (6.75,6.75);

\foreach \i in {.75, 1.5, 2.25, 3, 3.75, 4.5, 5.25} {
	\draw [thin, black] node at (\i+.375,6-\i+.375) {$*$};
	}

\foreach \i in {.75, 1.5, 2.25, 3} {
	\draw [thin, black] node at (\i+.375, 2.625) {$*$};
	}

\foreach \i in {0, 6} {
	\draw [thin, black] node at (\i+.375, 6-\i+.42) {$\ddots$};
	}

\draw [thin, black] node at (.375, 2.61) {$\cdots$};

\end{tikzpicture}
\end{center}

\noindent where * denotes a populated cell due to an initial condition

We can see that this is a valid initial conditions layout for this overlay because we can uniquely construct the array as we move around the overlay, filling in cells as shown below.

\begin{lemma}
If we have an overlay $B=(b_{i,j})$ with template $T=b_{0,0}+b_{1,0}Y+b_{1,1}XY$, where $b_{0,0}, b_{1,0}, b_{1,1} \ne 0$, then populated values of $A_{r,c}$ such that $A_{i,i}, i\in \Z$ and for a fixed row $k$, $A_{k,c}$ such that $c\in \Z, c<k$, then we can uniquely construct the array $A$ from this initial conditions layout.
\end{lemma}
\proof
Let $k$ be the integer such that $A_{k,c}$ such that $c<k$ is populated by the initial conditions. Since $T=b_{0,0}+b_{1,0}Y+b_{1,1}XY$, then for all $(r,c)\in \Z^2$, $$b_{0,0}A_{r,c}+b_{1,0}A_{r-1,c}+b_{1,1}A_{r-1, c-1}=0$$\\

\noindent Region 1:
First, we will show $A_{r,c}$ can be populated uniquely for $A_{r,c}$ such that $r \ge k, c\ge r$.
For the base case, we have that $A_{r,r}, r\in \Z$ is populated due to initial conditions.  In this case, $r=c$, or equivalently, $c-r=0$  Let $p \ge 0$, and suppose $A_{r,c}$ can be uniquely populated for $ 0 \le c-r \le p$.  We want to show this holds for $c-r=p+1$.  From the recurrence, we have that $b_{0,0}A_{r,c}+b_{1,0}A_{r-1,c}+b_{1,1}A_{r-1, c-1}=0$. Consider $c,r$ such that $c-r=p$, or $c=r+p$.  Then our recurrence becomes $$b_{0,0}A_{r,r+p}+b_{1,0}A_{r-1,r+p}+b_{1,1}A_{r-1, r+p-1}=0.$$  Notice that for the term $A_{r-1, r+p+1}$, subtracting the row index from the column index, we have $r+p-(r-1) = p+1$.  However, for the term $A_{r,r+p}$, subtracting the row index from the column index is $r+p-(r)=p$, so $A_{r,r+p}$ has already been populated by our assumption, and for the term $A_{r-1,r+p-1}$, subtracting the row index from the column index is $r+p-1 - (r-1)=p$, so $A_{r-1, r+p-1}$ has also already been populated by our assumption.  Since $b_{0,0}, b_{1,0}, b_{1,1} \ne 0$, then we have $$A_{r-1,r+p}=-\frac{1}{b_{1,0}} \left( b_{0,0}A_{r,r+p}+b_{1,1}A_{r-1, r+p-1} \right).$$  Thus $A_{r-1, r+p}$ is uniquely determined, thus $A_{r,c}$ is uniquely determined for $c-r=p+1$.\\

\noindent Region 2:
Next, we will show $A_{r,c}$ can be populated uniquely for $A_{r,c}$ such that $r \le k, c \le r$.  For the outer induction base case, we have that $A_{k,c}$ for $c \le k$ is populated due to initial conditions.  Furthermore, we also have that $A_{r,r}, r \in \Z$ is populated due to initial conditions.  Let $p \in \Z, p \ge 0$, and suppose, for the outer induction hypothesis, $A_{r,c}$ can be uniquely populated for $A_{k-p,c}, c\le k-p$.  We want to show we can uniquely populate $A_{k-(p+1),c}, c \le k-(p+1)$.   From the recurrence, we have that $b_{0,0}A_{r,c}+b_{1,0}A_{r-1,c}+b_{1,1}A_{r-1, c-1}=0$. For the inner induction base case, consider $r=k-p, c=k-(p+1)=k-p-1$.  Then the recurrence becomes $$b_{0,0}A_{k-p,k-p-1}+b_{1,0}A_{k-p-1,k-p-1}+b_{1,1}A_{k-p-1, k-p-2}=0.$$  Observe that $A_{k-p-1, k-p-1}$ is populated due to the initial conditions, and $A_{k-p, k-p-1}$ is populated due to the outer induction hypothesis since $c=k-p-1 < k-p$.  Thus since $b_{1,1} \ne 0$, we have $$A_{k-p-1, k-p-2}=-\frac{1}{b_{1,1}}\left(b_{0,0}A_{k-p, k-p-1}+b_{1,0}A_{k-p-1, k-p-1}\right),$$ so $A_{k-p-1, k-p-2}$ is uniquely populated by the recurrence.  Let $j\in \Z, j \ge 0$, and suppose, for the inner induction hypothesis, that $A_{r,c}$ can be uniquely populated for $A_{k-(p+1), k-(p+1)-j}$ i.e. $A_{k-p-1, k-p-j-1}$ (thus maintaining that the column index less than or equal to the row index).  We want to show $A_{k-(p+1), k-(p+1)-(j+1)}$ i.e. $A_{k-p-1, k-p-j-2}$ can be uniquely populated.  By considering $r=k-p, c=k-(p+1)-j$, the recurrence becomes $b_{0,0}A_{k-p,k-p-j-1}+b_{1,0}A_{k-p-1,k-p-j-1}+b_{1,1}A_{k-p-1, k-p-j-2}=0$ i.e. $$b_{0,0}A_{k-p,k-j-(p+1)}+b_{1,0}A_{k-(p+1),k-j-(p+1)}+b_{1,1}A_{k-(p+1), k-(p+1)-(j+1)}=0.$$  Since $A_{k-(p+1), k-j-(p+1)}$ is uniquely populated by the outer induction hypothesis and $A_{k-(p+1), k-j-(p+1)}$ is uniquely populated by the inner induction hypothesis.  Thus since $b_{1,1} \ne 0$, $$A_{k-(p+1), k-(p+1)-(j+1)} = -\frac{1}{b_{1,1}} \left(b_{0,0}A_{k-p,k-j-(p+1)}+b_{1,0}A_{k-(p+1),k-j-(p+1)}\right)$$ and thus $A_{k-(p+1), k-(p+1)-(j+1)}$ is uniquely determined by the recurrence.\\


\noindent Region 3:
Finally, we will show $A_{r,c}$ can be populated uniquely for $A_{r,c}$ such that $r \ge k, c \le r$.  For the outer induction base case, we have that $A_{k,c}$ for $c\le k$ is populated due to initial conditions.  Furthermore, we also have that $A_{r,r}, r \in \Z$ is populated due to initial conditions.  Let $p \in \Z, p \ge 0$, and suppose, for the outer induction hypothesis, $A_{r,c}$ can be uniquely populated for $A_{k+p,c}, c \le k+p$.  We want to show we can uniquely populate $A_{k+(p+1), c}, c \le k+(p+1)$.  From the recurrence, we have that $b_{0,0}A_{r,c}+b_{1,0}A_{r-1,c}+b_{1,1}A_{r-1, c-1}=0$.  For the inner induction base case, consider $r=k+p+1, c=k+p$.  Then the recurrence becomes $$b_{0,0}A_{k+p+1,k+p}+b_{1,0}A_{k+p,k+p}+b_{1,1}A_{k+p, k+p-1}=0.$$  Observe that $A_{k+p, k+p}$ is populated due to the initial conditions, and $A_{k+p, k+p-1}$ is populated due to the outer induction hypothesis since $c=k+p-1 < k+p$.  Thus since $b_{0,0} \ne 0$, we have $$A_{k+p+1, k+p}= -\frac{1}{b_{0,0}}\left(b_{1,0}A_{k+p,k+p}+b_{1,1}A_{k+p,k+p-1}\right),$$ so $A_{k+p+1, k+p}$ is uniquely populated by the recurrence.  Let $j \in \Z, j\ge 0$, and suppose, for the inner induction hypothesis, that $A_{r,c}$ can be uniquely populated for $A_{k+p+1, k+p+1-j}$ (thus maintaining that the column index less than or equal to the row index).  We want to show $A_{k+p+1, k+p+1 - (j+1)}$ i.e. $A_{k+p+1, k+p-j}$ can be uniquely populated.  By considering $r=k+p+1, c=k+p-j$, the recurrence becomes $b_{0,0}A_{k+p+1, k+p-j} + b_{1,0}A_{k+p,k+p-j} + b_{1,1}A_{k+p, k+p-j-1}=0$, i.e. $$b_{0,0}A_{k+p+1, k+p-j} + b_{1,0}A_{k+p,k+p-j} + b_{1,1}A_{k+p, k+p-(j+1)}=0.$$  Since $A_{k+p, k+p-j}$ is uniquely populated by the outer induction hypothesis since $c=k+p-j < k+p$ and $A_{k+p, k+p-(j+1)}$ is uniquely populated by the inner induction hypothesis. \\


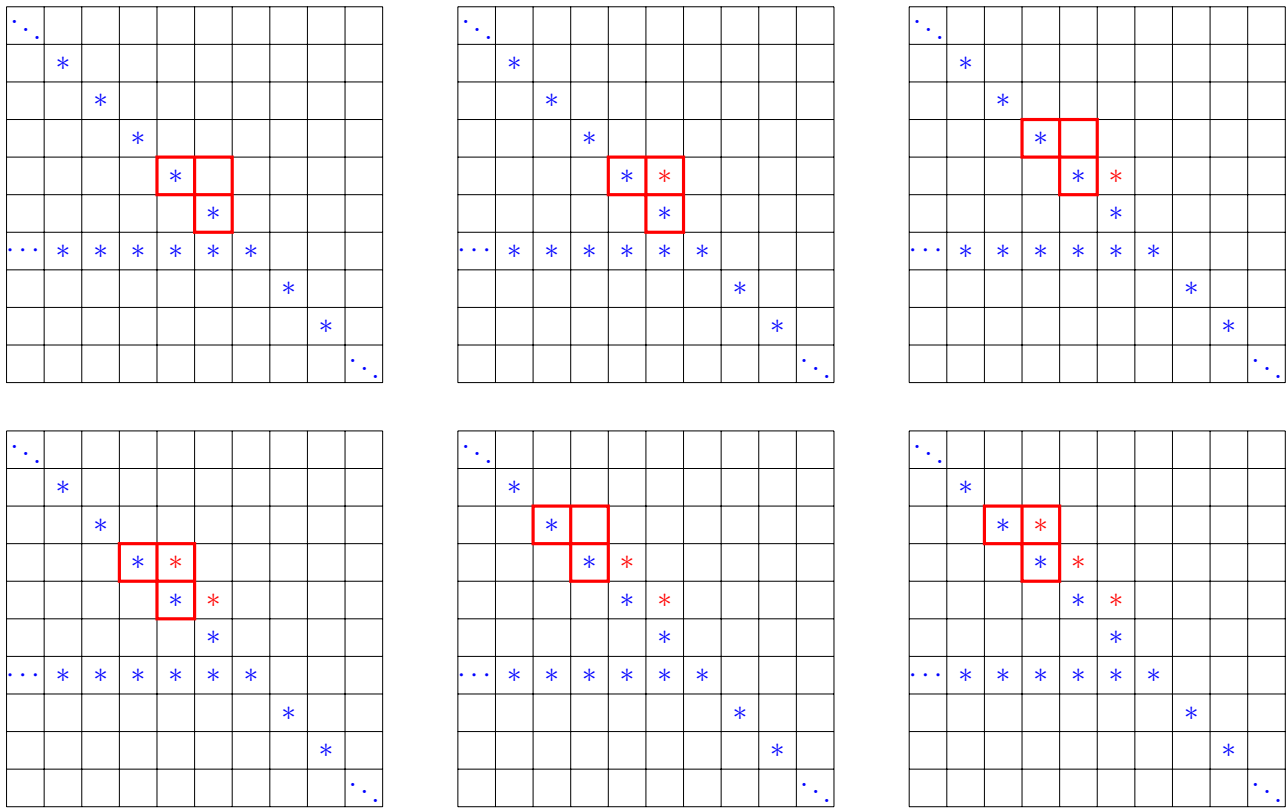
\begin{center}
\begin{tikzpicture}
{\small
\draw [step=.5, black, very thin] (0,0) grid (5,5);

\foreach \i in {.5, 1, 1.5, 2, 2.5, 3, 3.5, 4} {
	\draw [thin, blue] node at (\i+.25, 4.75-\i) {$*$};
	}

\foreach \i in {.5, 1, 1.5, 2, 2.5} {
	\draw [thin, blue] node at (\i+.25, 1.75) {$*$};
	}

\foreach \i in {0, 4.5} {
	\draw [thin, blue] node at (\i+.255, 4.5-\i+.3) {$\ddots$};
	}

\draw [thin, blue] node at (.25, 1.75) {$\cdots$};

\draw [very thick, red] (2,3)--(3,3)--(3,2)--(2.5,2)--(2.5,2.5)--(2, 2.5)--cycle;
\draw [very thick, red] (3,2.5)--(2.5,2.5)--(2.5,3);

\draw [step=.5, black, very thin] (5.999,0) grid (11,5);

\foreach \i in {.5, 1, 1.5, 2, 2.5, 3, 3.5, 4} {
	\draw [thin, blue] node at (\i+6.25, 4.75-\i) {$*$};
	}

\foreach \i in {.5, 1, 1.5, 2, 2.5} {
	\draw [thin, blue] node at (\i+6.25, 1.75) {$*$};
	}

\foreach \i in {0, 4.5} {
	\draw [thin, blue] node at (\i+6.255, 4.5-\i+.3) {$\ddots$};
	}

\draw [thin, blue] node at (6.25, 1.75) {$\cdots$};

\draw [very thick, red] (8,3)--(9,3)--(9,2)--(8.5,2)--(8.5,2.5)--(8, 2.5)--cycle;
\draw [very thick, red] (9,2.5)--(8.5,2.5)--(8.5,3);
\draw [thin, red] node at (8.75, 2.75) {$*$};

\draw [step=.5, black, very thin] (11.999,0) grid (17, 5);

\foreach \i in {.5, 1, 1.5, 2, 2.5, 3, 3.5, 4} {
	\draw [thin, blue] node at (\i+12.25, 4.75-\i) {$*$};
	}

\foreach \i in {.5, 1, 1.5, 2, 2.5} {
	\draw [thin, blue] node at (\i+12.25, 1.75) {$*$};
	}

\foreach \i in {0, 4.5} {
	\draw [thin, blue] node at (\i+12.255, 4.5-\i+.3) {$\ddots$};
	}

\draw [thin, blue] node at (12.25, 1.75) {$\cdots$};

\draw [very thick, red] (13.5,3.5)--(14.5,3.5)--(14.5,2.5)--(14,2.5)--(14,3)--(13.5, 3)--cycle;
\draw [very thick, red] (14.5,3)--(14,3)--(14,3.5);
\draw [thin, red] node at (14.75, 2.75) {$*$};

}
\end{tikzpicture}
\end{center}

\begin{center}
\begin{tikzpicture}
{\small
\draw [step=.5, black, very thin] (0,0) grid (5,5);

\foreach \i in {.5, 1, 1.5, 2, 2.5, 3, 3.5, 4} {
	\draw [thin, blue] node at (\i+.25, 4.75-\i) {$*$};
	}

\foreach \i in {.5, 1, 1.5, 2, 2.5} {
	\draw [thin, blue] node at (\i+.25, 1.75) {$*$};
	}

\foreach \i in {0, 4.5} {
	\draw [thin, blue] node at (\i+.255, 4.5-\i+.3) {$\ddots$};
	}

\draw [thin, blue] node at (.25, 1.75) {$\cdots$};

\draw [very thick, red] (1.5,3.5)--(2.5,3.5)--(2.5,2.5)--(2,2.5)--(2,3)--(1.5, 3)--cycle;
\draw [very thick, red] (2.5,3)--(2,3)--(2,3.5);
\draw [thin, red] node at (2.75, 2.75) {$*$};
\draw [thin, red] node at (2.25, 3.25) {$*$};

\draw [step=.5, black, very thin] (5.999,0) grid (11,5);

\foreach \i in {.5, 1, 1.5, 2, 2.5, 3, 3.5, 4} {
	\draw [thin, blue] node at (\i+6.25, 4.75-\i) {$*$};
	}

\foreach \i in {.5, 1, 1.5, 2, 2.5} {
	\draw [thin, blue] node at (\i+6.25, 1.75) {$*$};
	}

\foreach \i in {0, 4.5} {
	\draw [thin, blue] node at (\i+6.255, 4.5-\i+.3) {$\ddots$};
	}

\draw [thin, blue] node at (6.25, 1.75) {$\cdots$};

\draw [very thick, red] (7,4)--(8,4)--(8,3)--(7.5,3)--(7.5,3.5)--(7, 3.5)--cycle;
\draw [very thick, red] (8,3.5)--(7.5,3.5)--(7.5,4);
\draw [thin, red] node at (8.75, 2.75) {$*$};
\draw [thin, red] node at (8.25, 3.25) {$*$};

\draw [step=.5, black, very thin] (11.999,0) grid (17, 5);

\foreach \i in {.5, 1, 1.5, 2, 2.5, 3, 3.5, 4} {
	\draw [thin, blue] node at (\i+12.25, 4.75-\i) {$*$};
	}

\foreach \i in {.5, 1, 1.5, 2, 2.5} {
	\draw [thin, blue] node at (\i+12.25, 1.75) {$*$};
	}

\foreach \i in {0, 4.5} {
	\draw [thin, blue] node at (\i+12.255, 4.5-\i+.3) {$\ddots$};
	}

\draw [thin, blue] node at (12.25, 1.75) {$\cdots$};

\draw [very thick, red] (13,4)--(14,4)--(14,3)--(13.5,3)--(13.5,3.5)--(13, 3.5)--cycle;
\draw [very thick, red] (14,3.5)--(13.5,3.5)--(13.5,4);
\draw [thin, red] node at (14.75, 2.75) {$*$};
\draw [thin, red] node at (14.25, 3.25) {$*$};
\draw [thin, red] node at (13.75, 3.75) {$*$};

}
\end{tikzpicture}
\end{center}

\newpage


\begin{lemma}
Let $B=(b_{ij})$ be a $2 \times 2$ overlay and $E(i,j)$ be basis arrays built from the overlay applied to initial conditions defined as follows: for $k,l$ in the initial conditions layout $\gamma$,
\[ E(i,j)_{k,l} =
\begin{cases}
1 & k=i \text{ and } l=j \\
0 & otherwise.
\end{cases}
\]
For $(i,j) \in \gamma$, \\
\indent Case 1: If $j>0,$ then $E(i,j)_{k,l}=0$ for $l<j$. \\
\indent Case 2: If $j<0,$ then $E(i,j)_{k,l}=0$ for $l>j$. \\
\indent Case 3: If $i>0,$ then $E(i,j)_{k,l}=0$ for $k<i$. \\
\indent Case 4: If $i<0,$ then $E(i,j)_{k,l}=0$ for $k>i$. \\

\end{lemma}

\proof


Case 1: $j>0$.  Without loss of generality, take $\gamma$ to be the standard initial conditions layout.  Then  $E(0,j)_{\alpha, 0}=0$ for all $\alpha$ and $E(0,j)_{0,\beta}=0$ for all $\beta \ne j$ by the initial conditions. We will show $E(i,j)_{k,l}=0$ for $l<j$ by induction on $l$.  $l=0$ is given from the initial conditions, so assume that $E(0,j)_{k,l}=0$ for all $l$ such that $0\leq l \leq \lambda <j-1$.  We want to show this holds for $l=\lambda + 1$.

We know $E(0,j)_{0,\lambda+1}=0$ by the initial conditions, so assume by induction on $k>0$ that $E(0,j)_{k,\lambda+1}=0$.  We want to show that $E(0,j)_{k+1,\lambda+1} = 0.$  By our recurrence, we have

$$E(0,j)_{k+1,\lambda+1}=-\dfrac{1}{b_{00}}(b_{11}E(0,j)_{k,\lambda}+b_{10}E(0,j)_{k,\lambda+1}+b_{01}E(0,j)_{k+1,\lambda})$$
By our inductive hypothesis on $l$, $E(0,j)_{k,\lambda}=0$ and $E(0,j)_{k+1,\lambda}=0$, and by our inductive hypothesis on $k$, $E(0,j)_{k,\lambda+1}=0.$ Thus,
$$E(0,j)_{k+1,\lambda+1}=-\dfrac{1}{b_{00}}\left(b_{11}(0)+b_{10}(0)+b_{01}(0)\right)=0$$.

Similarly, assume by induction on $k<0$ that $E(0,j)_{k,\lambda+1}=0$.  We want to show that $E(0,j)_{k-1,\lambda+1} = 0.$  By our recurrence, we have
$$E(0,j)_{k-1,\lambda+1}=-\dfrac{1}{b_{10}}(b_{11}E(0,j)_{k-1,\lambda}+b_{00}E(0,j)_{k,\lambda+1}+b_{01}E(0,j)_{k,\lambda})$$
By our inductive hypothesis on $l$, $E(0,j)_{k,\lambda}=0$ and $E(0,j)_{k-1,\lambda}=0$, and by our inductive hypothesis on $k$, $E(0,j)_{k,\lambda+1}=0.$ 
Thus,
$$E(0,j)_{k-1,\lambda+1}=-\dfrac{1}{b_{1,0}}(b_{1,1}(0)+b_{0,0}(0)+b_{0,1}(0))=0$$.\\


Therefore for $j>0$, $E(0,j)_{k,l}=0$ for all $l<j$. \qed


\vspace{.5in}

\begin{theorem}
Let $\alpha$ be an initial conditions layout and let $$B=\sum_{(i,j)\in \alpha } a_{i,j} E(i,j) = \left(b_{k,l}\right),$$ where $a_{i,j}\in \F$.  There are only a finite number of $E(i,j)$ that contribute nonzero entries to each $b_{k,l}$.

\end{theorem}


\vspace{2in}
\begin{theorem}
Let $B=(b_{i,j})$ be an $(m+1) \times (n+1)$ overlay and $E(i,j)$ be basis arrays built from the overlay applied to initial conditions defined as follows: for $(k,l)$ in the initial conditions layout $\gamma$,
\[ E(i,j)_{k,l} =
\begin{cases}
1 & \text{if } k=i \text{ and } l=j \\
0 & otherwise.
\end{cases}
\]
For $(i,j) \in \gamma$, \\
\indent Case 1: If $j\ge m,$ then $E(i,j)_{k,l}=0$ for $l<j$. \\
\indent Case 2: If $j<0,$ then $E(i,j)_{k,l}=0$ for $l>j$. \\
\indent Case 3: If $i\ge n,$ then $E(i,j)_{k,l}=0$ for $k<i$. \\
\indent Case 4: If $i<0,$ then $E(i,j)_{k,l}=0$ for $k>i$. \\
\end{theorem}

\proof
Case 1: $j\ge m$.  Without loss of generality, take $\gamma$ to be the standard initial conditions layout.  Then  $E(0,j)_{\alpha, 0}=0$ for all $\alpha$ and $E(0,j)_{0,\beta}=0$ for all $\beta \ne j$ by the initial conditions. We will show $E(i,j)_{k,l}=0$ for $l<j$ by induction on $l$.  $l=0$ is given from the initial conditions, so assume that $E(0,j)_{k,l}=0$ for all $l$ such that $0\leq l \leq \lambda <j-1$.  We want to show this holds for $l=\lambda + 1$.

We know $E(0,j)_{0,\lambda+1}=0$ by the initial conditions, so assume by induction on $k>0$ that $E(0,j)_{k,\lambda+1}=0$.  We want to show that $E(0,j)_{k+1,\lambda+1} = 0.$  By our recurrence, we have

$$E(0,j)_{k+1,\lambda+1}=-\dfrac{1}{b_{00}}(b_{11}E(0,j)_{k,\lambda}+b_{10}E(0,j)_{k,\lambda+1}+b_{01}E(0,j)_{k+1,\lambda})$$
By our inductive hypothesis on $l$, $E(0,j)_{k,\lambda}=0$ and $E(0,j)_{k+1,\lambda}=0$, and by our inductive hypothesis on $k$, $E(0,j)_{k,\lambda+1}=0.$ Thus,
$$E(0,j)_{k+1,\lambda+1}=-\dfrac{1}{b_{00}}\left(b_{11}(0)+b_{10}(0)+b_{01}(0)\right)=0$$.

Similarly, assume by induction on $k<0$ that $E(0,j)_{k,\lambda+1}=0$.  We want to show that $E(0,j)_{k-1,\lambda+1} = 0.$  By our recurrence, we have
$$E(0,j)_{k-1,\lambda+1}=-\dfrac{1}{b_{10}}(b_{11}E(0,j)_{k-1,\lambda}+b_{00}E(0,j)_{k,\lambda+1}+b_{01}E(0,j)_{k,\lambda})$$
By our inductive hypothesis on $l$, $E(0,j)_{k,\lambda}=0$ and $E(0,j)_{k-1,\lambda}=0$, and by our inductive hypothesis on $k$, $E(0,j)_{k,\lambda+1}=0.$ 
Thus,
$$E(0,j)_{k-1,\lambda+1}=-\dfrac{1}{b_{1,0}}(b_{1,1}(0)+b_{0,0}(0)+b_{0,1}(0))=0$$.\\

Therefore for $j>0$, $E(0,j)_{k,l}=0$ for all $l<j$. \qed


\begin{lemma}
\[A=\sum_{(i,j) \in \gamma} d_{i,j}E(i,j)\]
where $d_{i,j}$ is the initial condition of the array $A$ in the $(i,j)$ position.
Assuming that for any position $(r,c)$, there are finitely many $E(i,j)_{r,c}$ that contribute nonzero values.
\end{lemma}

\proof 

Consider $$\left( \sum_{(i,j)\in \gamma} d_{i,j}E(i,j) \right)_{r,c}$$



 \[ \left(B \sum_{(i,j)\in \gamma} d_{i,j}E(i,j) \right)_{r,c} = \sum_{(i,j)\in \gamma} d_{i,j}B \left(E(i,j)\right)_{r,c} \]

Since $E(i,j)$ is annihilated by $B$, then 
\[ \left( \sum_{(i,j)\in \gamma} d_{i,j}B \left(E(i,j)\right) \right)_{r,c} = 0\].

\bibliography{Templates_Arrays_and_Overlays} 
\bibliographystyle{plain}

\end{document}